\documentclass[a4paper, 11pt]{article}
\usepackage{amsmath, amsfonts,amssymb}
\usepackage[english]{babel}
\begin {document}

\begin{center}
\textbf {Geometrical properties of the space $P_f(X)$ \\ of probability measures} \\

\medskip
\textbf {Zaitov A.\,A.}\\
\smallskip
{Tashkent institute of architecture and civil engineering}\\
{adilbek\_zaitov@mail.ru},\\

\end{center}

\begin{abstract}
In this paper we prove that for a compact space $X$ inclusion $P_{f}(X)\in ANR$ holds if and only if $X\in ANR$. Further, it is shown that the functor $P_{f}$ preserves property of a compact to be  $Q$-manifold or a Hilbert cube, properties of maps fibres to be $ANR$-compact, $Q$-manifold, Hilbert cube (the finite of Hilbert cube).\\

\end{abstract}

Creation of the general theory of infinite-dimensional manifolds has increased interest in infinite-dimensional objects existing in `nature'. The functor $P$,  transforming arbitrary compact to convex subsets of local convex spaces, supplies such objects. Among other covariant functors the functor $P$ distinguishes with the fact that investigation of probability measures is conducted on a joint at least of three mathematical disciplines: topology, functional analysis and probability theory. What explains a variety of the applied methods and great opportunities for applications of the received results.

In work [1] the subfunctor $P_{f}$ of the functor $P$ of probability measures was entered. For a compact $X$ the set  $P_{f}(X)$ consists of probability measures with finite support, and if the support of a measure $\mu$ consists of $n$ points $x_{1}, x_{2}..., x_{n}$, the barycenter mass of one of these points isn't less than $1-\frac{1}{n+1}$. This functor is interesting that it is a functor with the finite support, and has no final degree. Functor $P_{f}: Comp\rightarrow  Comp$ is a normal subfunctor of functor $P$ of probability measures. Earlier in works [2-5] functor $P_{f}$ it has been investigated. In the present work we state strict proofs of the received results.

From definition of elements of the space $P_{f}(X)$ follows that a set $\delta(X)$ of Dirac measures lies in $P_{f}(X)$.

Let $X$ and $Y$ be two compacts lying in spaces $M$ and $N$  respectively, where  $M,N\in AR$. A sequence of maps $f_{k}:M\rightarrow N$, $k=1,2,...,$ is called to be fundamental sequence from $X$ into $Y$, if for each neighbourhood  $V$ of the compact $Y$  (in $N$) there is such neighbourhood $U$ of the compact $X$   (in $M$), that $f_{k}|_{U}=f_{k+1}|_{U}$  at $V$ almost for all $k$. It means that there is such homotopy $f_{k}:U\times[0,1]\rightarrow V$, that $f_{k}(x,0)=f_{k}(x)$ and $f_{k+1}(x,1)=f_{k+1}(x)$  for all $x\in U$. We will denote this fundamental sequence through $\left\{f_{k},X,Y\right\}$  or shortly through $f$, also we will write $f:X\rightarrow Y$.

A fundamental sequence $f=\left\{f_{k},X,Y\right\}$ is generated by map $f:X\rightarrow Y$ if  $f_{k}(x)=f(x)$  for all $x\in X$ and for all $k=1,2,...$.

Spaces $X$ and $Y$ are fundamentally equivalent if there are such two fundamental sequences $f:X\rightarrow Y$ and $g:Y\rightarrow X$ that $gf=id_{X}$ and $fg=id_{Y}$.

The relation of fundamental equivalence is the equivalence relation therefore the class of all spaces decomposes to in pairwise disjoint classes of spaces which are called shape. So two spaces belong to the same shape if and only if when they are fundamentally equivalent. A shape containing a space $X$ is called [6] a shape of the space $X$ and denote by $Sh(X)$. It is known that for two neighbourhood retracts $A$ and $B$ the equality $Sh(A)=Sh(B)$ is true if and only if when they are homotopically equivalent.

Let's remind that maps  $f:X\rightarrow Y$ is called cellularity-similar (briefly $CE$) [7], if any compact $A\subset Y$  preimage $f^{-1}(A)$  is a compact and for each point $y\in Y$  preimage $f^{-1}(y)$  has shape points (i. e. a preimage $f^{-1}(y)$ is homotopically equivalent to a point).

If $r:X\rightarrow F$ is a retraction and there also exists such a homotopy $h:X\times[0,1]\rightarrow F$ that $h(x,0)=x$, $h(x,1)=r(x)$ for all $x\in X$ then $r$ is deformation retraction, and $F$ is deformation retract of the space $X$. Deformation retraction  $r:X\rightarrow F$  is strongly deformation retraction if for a homotopy $h:X\times[0,1]\rightarrow F$ we have $h(x,t)=x$ for all $x\in F$  and all $t\in[0,1]$  [6].

$\mu=\bigoplus\limits_{i=1}^{n}\lambda_{i}\odot\delta_{x_{i}}\in I_{f}(X)$. Let $\lambda_{i_{0}}=0>-ln(n+1)$. Measure $\mu$  puts in correspondence a point $\delta_{x_{i_{0}}}$ of a compact $\delta(X)$. The obtained correspondence $I_{f}(X)\rightarrow\delta(X)$  denote by $r_{\delta(X)}^{I(X)}$.

Take arbitrary measure $\mu\in P_{f}(X)$, $\mu=\sum\limits_{i=1}^{n}m_{i}\delta_{x_{i}}$. Let $m_{i_{0}}\geq 1-\frac{1}{n+1}$. Measure $\mu$ assign a point $\delta_{x_{i_{0}}}$ of the compact $\delta(X)$. The obtained correspondence $r_{\delta(X)}^{P_f(X)}: P_{f}(X)\rightarrow\delta(X)$ denote by $r_{\delta}^{f}$.

{\bf Theorem 1.} {\it For any compact $X$ a map $r_{\delta}^{f}:P_{f}(X)\rightarrow \delta(X)$ is continuous, open, cellularity-similar (all fibres are collapsible) retraction.}

{\bf Proof.} By construction we have $\left(r_{\delta}^{f}\right)^{-1}(\delta_z)\cap \left(r_{\delta}^{f}\right)^{-1}(\delta_y)=\varnothing$ for every pair $y,\ z\in X$, $y\neq z$. Therefore the map $r_{\delta}^{f}:P_{f}(X)\rightarrow \delta(X)$ is defined correctly. It is clear that $r_{\delta}^{f}(\delta_{x})=\delta_{x}$ for each  $x\in X$, i. e. every point of the space $\delta(X)$ is fixed-point according to the map  $r_{\delta}^{f}:P_{f}(X)\rightarrow \delta(X)$. So, we establish that $r_{\delta}^{f}$ is retraction.

It is clear, that for each point $x\in X$ the fibre $\left(r_{\delta}^{f}\right)^{-1}(\delta_{x})$ is compact. On the other hand for each point $\mu\in \left(r_{\delta}^{f}\right)^{-1}(\delta_{x})$ an interval $[\mu,\delta_{x}]=\{\alpha\delta_{x}+(1-\alpha)\mu: 0\leq\alpha\leq 1\}$ lies in the fibre $\left(r_{\delta}^{f}\right)^{-1}(\delta_{x})$.

Fix a fibre $\left(r_{\delta}^{f}\right)^{-1}(\delta_{x})$ and define a map
$h:\left(r_{\delta}^{f}\right)^{-1}(\delta_{x})\times[0,1]\rightarrow \left(r_{\delta}^{f}\right)^{-1}(\delta_{x})$  by the rule
$$h(\mu,t)=(1-t)\delta_{x}+t\mu,$$
where $\mu=\sum\limits_{i=1}^{n}m_i\delta_i\in \left(r_{\delta}^{f}\right)^{-1}(\delta_{x})$ и $t\in [0,\ 1].$

$h$  ,  $h_1=id_{\left(r_{\delta(X)}^{I(X)}\right)^{-1}(\delta_x)}$ $h_0=\left(r_{\delta(X)}^{I(X)}\right)^{-1}(\delta_x)\rightarrow\{\delta_x\}$.

It is easy to see that $h$ is a homotopy, connecting stationary map $h_0:\left(r_{\delta}^{f}\right)^{-1}(\delta_{x})\rightarrow\{\delta_x\}$ and identity mapping  $h_1=id_{\left(r_{\delta}^{f}\right)^{-1}(\delta_{x})}$.

So each pre-image $\left(r_{\delta}^{f}\right)^{-1}(\delta_{x})$, $x\in X$, has shape of a point, i. e. a retraction $r_{\delta}^{f}$  is cellularity-similar. Moreover, these pre-images subtend to a point.

We will show a map $r_{\delta}^{f}:P_{f}(X)\rightarrow\delta(X)$ is continuous. Note that at variation of measures $\mu_{0}\in P(X)$, of finite collections $\{U_{1},...,U_{n}\}$ co-zero (respectively, open) sets $U_{i}$ of the space (respectively, of the compact) $X$ and of numbers $\varepsilon >0$, a family of sets of the view
$$\langle\mu_{0};U_{1},...,U_{n};\varepsilon\rangle=\{\mu\in P(X):\mu(U_{i})-\mu_{0}(U_{i})>-\varepsilon, \ \ i=1,...,n\}$$ forms [8] a base of the weak convergence topology of the space $P(X)$.

It is clear that $r_{\delta}^{f}(\langle\delta_{x};U_{1},...,U_{n};\varepsilon\rangle\cap\delta(X))= \langle\delta_{x};U_{1},...,U_{n};\varepsilon\rangle\cap\delta(X)$, i. e. $r_{\delta}^{f}|_{\delta(X)}$ is open map. That is why, further we will consider measures which support consists not less than two points and also without losing generality, we will suppose $\varepsilon <\frac{1}{3}$.

Let $\mu_{0}=\sum\limits_{i=1}^{k}m_{0i}\delta_{x_{i}}\in P_{f}(X), \ \ m_{0i_{0}}\geq 1-\frac{1}{k+1}, \ \ r_{\delta}^{f}(\mu_{0})=\delta_{x_{i_{0}}}$, and $V$ be an open set such that $r_{\delta}^{f}(\mu_{0})\in\langle V\rangle=\{\delta_{x}:\delta_{x}(V)>0\}$. Consider a (subbase's) neighbourhood 
$$\langle\mu_{0};V;\varepsilon\rangle=\{\mu\in P(X):\mu(V)-\mu_{0}(V)>-\varepsilon\}\cap P_{f}(X)$$
of the measure $\mu_{0}$ where $0<\varepsilon<\frac{1}{3}$.

Let $\mu=\sum\limits_{i=1}^{l}m_{i}\delta_{x_{i}}\in\langle\mu_{0};V;\varepsilon\rangle$. On construction 
$m_{i'}\geq1-\frac{1}{l+1}$ for some $i'$. It is evident that $\mu(V)\geq m_{i'}\geq1-\frac{1}{l+1}$ if and only if $x_{i'}\in V$. At $x_{i'}\notin V$ one has $\mu(V)\leq\frac{1}{l+1}$, and hence, $\mu(V)-\mu_{0}(V)\leq\frac{1}{l+1}-\frac{k}{k+1}\leq-\frac{1}{3}<-\varepsilon$, contrary to $\mu\in\langle\mu_{0};V;\varepsilon\rangle$. Whence it follows that $x_{i'}\in V$ i. e. $r_{\delta}^{f}(\mu)=\delta_{x_{i'}}\in\langle V\rangle$. In other words, we establish that $r_{\delta}^{f}(\langle\mu_{0};V;\varepsilon\rangle)\subset\langle V\rangle$. Thus the map $r_{\delta}^{f}$  is continuous.

Further, we will show that $r_{\delta}^{f}(\langle\mu_{0};V;\varepsilon\rangle)=\langle V\rangle=\{\delta_{x}:\delta_{x}(V)>0\}$. For arbitrary point $x\in V$ we will construct a measure $\mu_{x}$ as follows: $\mu_{x}=m_{0i_{0}}\delta_{x}+\sum\limits_{\substack{i=1,\\i\neq i_{0}}}^{k}m_{0i}\delta_{x_{i}}$. Then 
$\mu_{0}(V)-\mu_{x}(V)=0$, and hence $\mu_{x}\in\langle\mu_{0};V;\varepsilon\rangle$ for any $\varepsilon>0$. In other words for each $\delta_{x}\in\langle V\rangle$ we have established that $(r_{\delta}^{f})^{-1}(\delta_{x})\cap\langle\mu_{0};V;\varepsilon\rangle\neq\varnothing$. Whence it follows that $r_{\delta}^{f}(\langle\mu_{0};V;\varepsilon\rangle)=\langle V\rangle$. Thus, the image $r_{\delta}^{f}(\langle\mu_{0};V;\varepsilon\rangle)$ of open in $P_{f}(X)$ set $\langle\mu_{0};V;\varepsilon\rangle$ is open in $\delta(X)$. Theorem 1 is proved.

{\bf Proposition 1.} {\it For a compact $X$ subspace $\delta(X)$ is strongly deformation retract of the compact $P_{f}(X)$.}

{\bf Proof.} Consider a map $h:P_{f}(X)\times[0,1]\rightarrow P_{f}(X)$, defined by a formula $$h(\mu,t)=h_{t}(\mu)=(1-t)\cdot\mu+t\cdot r_{\delta}^{f}(\mu), \ \ \ \ (\mu,t)\in P_{f}(X)\times[0,1].$$
It easy to check that the map $h$ is defined correctly. Moreover $h_{0}=id_{P_{f}(X)}$ and $h_{1}=r_{\delta}^{f}$, i. e. $h$ is a homotopy, connecting the maps $id_{P_{f}(X)}$ and $r_{\delta}^{f}$.
Else we have  
$$h(\delta_{x},t)=(1-t)\cdot\delta_{x}+t\cdot r_{\delta}^{f}(\delta_{x})=(1-t)\cdot\delta_{x}+t\cdot \delta_{x}=\delta_{x},$$
i. e. $h_{t}(\delta_{x})=\delta_{x}$ for all $\delta_{x}\in\delta(X)$ and $t\in[0,1]$. Thus, $\delta(X)$  is strongly deformation retract of the compact $P_{f}(X)$. Proposition 1 is proved.

{\bf Proposition 2.} {\it For any finite compact $X$ the set $P_{f}(X)$ is neighbourhood retract of the compact  $P(X)$.}

{\bf Proof.} For each Dirac measure $\delta_{x}$, $x\in X$, we will construct neighbourhood $$\langle\delta_{x};\frac{1}{2}\rangle=\left\{\mu\in P(X):\mu(U)>\frac{1}{2} \mbox{ for each open set } U \mbox{ containing } x \right\},$$
and consider a set $\bigcup\limits_{x\in X}\langle\delta_{x};\frac{1}{2}\rangle$ which is open in $P(X)$. It is obviously that $P_{f}(X)\subset\bigcup\limits_{x\in X}\langle\delta_{x};\frac{1}{2}\rangle$.

Besides, $\langle\delta_{x};\frac{1}{2}\rangle\cap\langle\delta_{y};\frac{1}{2}\rangle=\varnothing$ at $x\neq y$. Really, let $x\neq y$ and $\mu\in\langle\delta_{x};\frac{1}{2}\rangle\cap\langle\delta_{y};\frac{1}{2}\rangle$. Then $\mu(U)>\frac{1}{2}$ and $\mu(V)>\frac{1}{2}$ for all open sets $U\ni x$ and $V\ni y$. In particular, for any pair of disjoint open sets $U$ and $V$ containing $x$ and $y$, respectively, we have $\mu(U)>\frac{1}{2}$ and $\mu(V)>\frac{1}{2}$. Hence, $\mu(X)\geq\mu(U\cup V)=\mu(U)+\mu(V)>\frac{1}{2}+\frac{1}{2}=1$ contrary to $\mu\in P(X)$. Thus $\langle\delta_{x};\frac{1}{2}\rangle\cap\langle\delta_{y};\frac{1}{2}\rangle=\varnothing$ at $x\neq y$.

We will show that $P_{f}(X)$ is retract of the open set $\bigcup\limits_{x\in X}\langle\delta_{x};\frac{1}{2}\rangle$. If $\nu\in\bigcup\limits_{x\in X}\langle\delta_{x};\frac{1}{2}\rangle$, where $\nu=\sum\limits_{i=1}^{n}n_{i}\delta_{x_{i}}$,  $\sum\limits_{i=1}^{n}n_{i}=1$,  $n_{i}\geq0$,  $i=1,...,n$, then it is obvious that $n_{i_{0}}>\frac{1}{2}$ for some (a unique) $i_{0}$, and hence $r_{\delta}^{f}(\nu)=\delta_{x_{i_{0}}}$.

Determine a map $r:\bigcup\limits_{x\in X}\langle\delta_{x};U;\frac{1}{2}\rangle\rightarrow P_{f}(X)$ by the rule
\[
r(\nu) =
\begin{cases}
\frac{n}{n+1}\cdot\delta_{x_{i_{0}}}+\sum\limits_{\substack{i=1,\\ i\neq i_{0}}}^{n}\frac{n_{i}}{(n+1)(1-n_{i_{0}})}\delta_{x_{i}}, & \text{at $n_{i_{0}}<\frac{n}{n+1}$;} \\
\sum\limits_{i=1}^{n}n_{i}\delta_{x_{i}}, & \text{at $n_{i_{0}}\geqslant \frac{n}{n+1}$.}
\end{cases}
\]
The map $r$ is correctly defined. It is continuous. Besides $r(\mu)=\mu$ for any measure $\mu\in P_{f}(X)$. Thus, the map $r$ is retraction. Proposition 2 proved.

Remind [7] that a set $A\subset X$  is called to be collapsible by space $X$ to a set $B\subset X$ if embedding map $i_{A}:A\rightarrow X$ is homotopic to some map $f:A\rightarrow X$ such that $f(A)\subset B$. If $B$  consists of one point, they say that $A$ is collapsible by $X$.

It is clear, if there is a homotopy $h:A\times I\rightarrow A$, such that $h(y,0)=i_{A}$, and $h(y,1)=\{\text {a point}\}$ then $A$ collapsible by $X$.

A space $X$  is called to be locally collapsible to a point $x_{0}\in X$ if any neighbourhood $U$ of the point $x_{0}$ contains a neighbourhood $U_{0}$ such, that collapsible to a point by $U$. A space $X$ is called to be locally collapsible if it is locally collapsible to each its point.

{\bf Theorem 2.} {\it Functor $P_{f}$ preserves collapsibility of compacts, i. e. if $X$ is a collapsible compact then $ P_{f}(X)$ is also collapsible compact.}

{\bf Proof.} We will show more: the functor $P_{f}$ preserves homotopy of maps. Let $h_{0},h_{1}:X\rightarrow Y$ be homotopic maps, $h:X\times[0,1]\rightarrow Y$ be a homotopy connecting maps $h_{0}, \ \ h_{1}$, i. e. $h(x,0)=h_{0}(x)$,  $h(x,1)=h_{1}(x)$. An embedding $i_{t_{0}}:X\times\{t_{0}\}\rightarrow X\times I$, defined by equality $i_{t_{0}}(x)=(x,t_{0})$,  $x\in X$, generates an embedding $P_{f}(i_{t_{0}}):P_{f}(X\times\{t_{0}\})\rightarrow P_{f}(X\times I)$.
But for every $t_{0}\in[0,1]$ a space $P_{f}(X\times\{t_{0}\})$ naturally homeomorphic to $P_{f}(X)\times\{t_{0}\}$. This homeomorphism one can carry out, as it is easy to see, by means of correspondence $\mu_{t_{0}}\leftrightarrow(\mu,t_{0})$, where $\mu_{t_{0}}=\sum\limits_{i=1}^{n}m_{i}\delta_{(x_{i},t_{0})}\in P_{f}(X\times\{t_{0}\})$ and $\mu=\sum\limits_{i=1}^{n}m_{i}\delta_{x_{i}}\in P_{f}(X)$.

Now, determine a map $P_{f}(h):P_{f}(X)\times[0,1]\rightarrow P_{f}(Y)$ by equality
$$P_{f}(h)\left(\sum\limits_{i=1}^{n}m_{i}\delta_{x_{i}},t\right)=\sum\limits_{i=1}^{n}m_{i}\delta_{h(x_{i},t)}.$$
We have 
$$P_{f}(h)\left(\sum_{i=1}^{n}m_{i}\delta_{x_{i}},0\right)= \sum_{i=1}^{n}m_{i}\delta_{h(x_{i},0)}=\sum_{i=1}^{n}m_{i}\delta_{h_{0}(x_{i})}=
P_{f}(h_{0})\left(\sum_{i=1}^{n}m_{i}\delta_{x_{i}}\right),$$ $$P_{f}(h)\left(\sum_{i=1}^{n}m_{i}\delta_{x_{i}},1\right)= \sum_{i=1}^{n}m_{i}\delta_{h(x_{i},1)}=\sum_{i=1}^{n}m_{i}\delta_{h_{1}(x_{i})}= P_{f}(h_{1})\left(\sum_{i=1}^{n}m_{i}\delta_{x_{i}}\right),$$
i. e. $P_{f}(h)(\mu,0)=P_{f}(h_{0})(\mu)$ and $P_{f}(h)(\mu,1)=P_{f}(h_{1})(\mu)$ for any $\mu\in P_{f}(X)$. In other words $P_{f}(h)$ is homotopy connecting maps $P_{f}(h_{0})$ and $P_{f}(h_{1})$. Thus, functor $P_{f}$ preserves homotopy of maps. Theorem 2 proved.

{\bf Lemma 1.} {\it For a compact $X$ the set $P_{f}(X)$ is the neighbourhood retract of the space  $P_{\omega}(X)$.}

{\bf Proof.} Construct the following subsets$$P_{\frac{\delta}{2}}(X)=\left\{\frac{1}{2}\delta_{x}+\frac{1}{2}\delta_{y}:x,y\in X, x\neq y\right\}\subset P_{2}(X),$$
$$P_{\frac{\omega}{2}}(X)=\left\{\sum_{i=1}^{n}m_{i}\mu_{i}:\mu_{i}\in P_{\frac{\delta}{2}}(X);\sum_{i=1}^{n}m_{i}=1,  m_{i}\geq0\right\}\subset P_{\omega}(X).$$

{\bf Claim 1.} {\it Для компакта $X$ the set $P_{\frac{\delta}{2}}(X)$ is compact.}

{\bf Proof.} Let $\xi=\theta\delta_{x}+(1-\theta)\delta_{y}\in\left [P_{\frac{\delta}{2}}(X)\right]_{P_{2}(X)}$,  $0<\theta<1$. Then there exists a net $\left\{\xi_{\alpha}= \frac{1}{2}\delta_{x_{\alpha}} +\frac{1}{2}\delta_{y_{\alpha}}\right\}\subset P_{\frac{\delta}{2}}(X)$ converging to $\xi$  according to weak convergence topology , i. e. $\xi_{\alpha}(\varphi) \xrightarrow [\  \alpha\  ] {} \xi(\varphi)$ for all function $\varphi\in C(X)$. Since $x\neq y$ and $X$ is Hausdorff space, there exist neighbourhoods  $Ox$ and $Oy$ such that $Ox\cap Oy=\varnothing$. That is why without losing generality one can suppose that $\left\{x_{\alpha}\right\}\subset Ox$  and  $\left\{y_{\alpha}\right\}\subset Oy$. Consider a function $\varphi\in C(X)$ such that $\varphi(x)=\varphi(x_{\alpha})=1$ and $\varphi(y)=\varphi(y_{\alpha})=0$ at all $\alpha$. For each $\varepsilon>0$ there exists $\alpha(\varepsilon)$ such that $\left|(\xi-\xi_{\alpha})(\varphi)\right|<\varepsilon$ at all $\alpha\succ\alpha(\varepsilon)$. It follows that $$\left|\left(\theta\delta_{x}+(1-\theta)\delta_{y}-\frac{1}{2}\delta_{x_{\alpha}}- \frac{1}{2}\delta_{y_{\alpha}}\right)(\varphi)\right|= \left|\theta-\frac{1}{2}\right|<\varepsilon.$$
Since $\varepsilon$ is arbitrary and the right part of the last equality does not from $\alpha$ we obtain equality $\theta=\frac{1}{2}$. Hence $\xi=\frac{1}{2}\delta_{x}+\frac{1}{2}\delta_{y}\in P_{\frac{\delta}{2}}(X)$. It follows that $\left[P_{\frac{\delta}{2}}(X)\right]_{P_{2}(X)}\subset P_{\frac{\delta}{2}}(X)$, i. e. $P_{\frac{\delta}{2}}(X)$ is closed in $P_{2}(X)=\{\mu\in P(X):\ |\mbox{supp}\mu|\leq 2\}$. Then $P_{\frac{\delta}{2}}(X)$ is compact by virtue of compactness of $P_{2}(X)$. Claim 1 is proved.

{\bf Claim 2.} {\it The set $P_{\frac{\omega}{2}}(X)$ is closed in $P_{\omega}(X)$.}

{\bf Proof.} Let $\mu=\sum_{k=1}^sm_{k}\delta_{x_{k}}\in\left[P_{\frac{\omega}{2}}(X)\right]_{P_{\omega}(X)}$,  and $\left\{\mu_{\alpha}\right\}\subset P_{\frac{\omega}{2}}(X)$ be a net converging to $\mu$ according to weak convergence topology, i. e. $\mu_{\alpha}(\varphi)\rightarrow\mu(\varphi)$ for all function $\varphi\in C(X)$. Let $\mu_{\alpha}=\sum\limits_{j=1}^{s_{\alpha}}m_{kl j}^{\alpha}\left(\frac{1}{2}\delta_{x_{kj}^{\alpha}}+\frac{1}{2}\delta_{x_{lj}^{\alpha}}\right)$, where $\{x_{kj}^{\alpha}\}\subset X$ is a net convergence to $x_{k}$ in topology of space $X$; here $k,l\in\{1,...,s\}$, $m_{kl j}^{\alpha}=m_{lk j}^{\alpha}$, $x_{kj}^{\alpha}\neq x_{lj}^{\alpha}$ at $k\neq l$ (note that to the same point $x_{k}$ can converge several nets, for examples, $x_{k j_{1}}^{\alpha}{\longrightarrow}_{\alpha}x_{k}$,  $x_{k j_{2}}^{\alpha}{\longrightarrow}_{\alpha}x_{k}$, and so on; and also, can it happens that a net is stationary). Since $x_{k}\neq x_{l}$ at $k\neq l$ and $X$ is Hausdorff space then there exist neighbourhoods  $Ox_{k}$ of points $x_{k}$,  $k=1,...,s$, such that $Ox_{k}\cap Ox_{l}=\varnothing$ at $k\neq l$. Therefore without losing generality one can take that $\left\{x_{k{j}}^{\alpha}\right\}_{\alpha}\subset Ox_{k}$,  $k=1,...,s$,  $j\in\left\{1,...,s_{\alpha}\right\}$. We have $\mu_{\alpha}= \sum\limits_{j=1}^{s_{\alpha}}m_{klj}^{\alpha}\left(\frac{1}{2}\delta_{x_{kj}^{\alpha}}+ \frac{1}{2}\delta_{x_{lj}^{\alpha}}\right)= \sum\limits_{k=1}^{s}\frac{\sum\limits_{l}m_{klj}^{\alpha}}{2}\delta_{x_{kj}^{\alpha}}$. Consider a function $\varphi\in C(X)$, such that $\varphi(x_{k})=\varphi\left(x_{kj}^{\alpha}\right)= sign\left(m_{k}-\frac{\sum\limits_{j,l}m_{klj}^{\alpha}}{2}\right)$ at all $\alpha$. Here $k=1,...,s$,  $j\in\{1,...,s_{\alpha}\}$. For each $\varepsilon>0$ there exists $\alpha(\varepsilon)$ such that 
$$\left|(\mu-\mu_{\alpha})(\varphi)\right|= \left|\left(\sum\limits_{k=1}^{s}m_{k}\delta_{x_{k}}- \sum_{k=1}^{s}\frac{\sum\limits_{j,l}m_{klj}^{\alpha}}{2}\delta_{x_{kj}^{\alpha}}\right)(\varphi)\right|=$$
$$\left|\sum\limits_{k=1}^{s}m_{k}\varphi(x_{k})- \sum_{k=1}^{s}\frac{\sum\limits_{j,l}m_{klj}^{\alpha}}{2}\varphi\left(x_{kj}^{\alpha}\right)\right|= \left|\sum\limits_{k=1}^{s}\left(m_{k}- \frac{\sum\limits_{j,l}m_{klj}^{\alpha}}{2}\right)\varphi\left(x_{k}\right)\right|$$
$$=\left|\sum\limits_{k=1}^{s}\left|m_{k}-\frac{\sum\limits_{j,l}m_{klj}^{\alpha}}{2}\right|\right|= \sum\limits_{k=1}^{s}\left|m_{k}-\frac{\sum_{j,l}m_{klj}^{\alpha}}{2}\right|<\varepsilon$$
at  $\alpha\succ\alpha(\varepsilon)$. That means,  $m_{k}=\lim\limits_{\alpha}\frac{\sum\limits_{j,l}m_{klj}^{\alpha}}{2}$. Suppose   $\lim\limits_{\alpha}m_{klj}^{\alpha}=m_{klj}$,  $m_{kl}=\sum\limits_{j}m_{klj}$. Then   $m_{k}=\frac{\sum\limits_{l}m_{kl}}{2}$. Hence,
$$\mu=\sum\limits_{k=1}^{s}m_{k}\delta_{x_{k}}= \sum\limits_{k=1}^{s}\frac{\sum\limits_{l}m_{kl}}{2}\delta_{x_{k}}= \sum\limits_{k,l}m_{kl}\left(\frac{1}{2}\delta_{x_{k}}+\frac{1}{2}\delta_{x_{l}}\right),$$
i. e.  $\mu\in P_{\frac{\omega}{2}}(X)$. Thus,  $P_{\frac{\omega}{2}}(X)$ is closed in $P_{\omega}(X)$. Claim 2 is proved.

{\bf Claim 3.} {\it Let $\mu=\sum\limits_{i=1}^{s}m_{i}\delta_{x_{i}}\in P_{\omega}(X)$. Inclusion $\mu\in P_{\frac{\omega}{2}}(X)$ is true if and only if $m_{i}\leq\frac{1}{2}$ at all $i=1,...,s$.}

{\bf Proof.} By virue of construction every measure $\mu\in P_{\frac{\omega}{2}}(X)$ represents in the view $\mu=\sum\limits_{1\leq i,j\leq s}m_{ij}\left(\frac{1}{2}\delta_{x_{i}}+ \frac{1}{2}\delta_{x_{j}}\right)$. Since $m_{ij}=m_{ji}$, then we have $\mu=\sum\limits_{i=1}^{s}\frac{\sum\limits_{j}m_{ij}}{2}\delta_{x_{i}}$. Denote  $m_{i}=\frac{1}{2}\sum\limits_{j}m_{ij}$, and by virtue of $\sum\limits_{j}m_{ij}\leq \sum\limits_{1\leq i, j\leq s}m_{ij}=1$, we obtain $m_{i}\leq \frac{1}{2}$ for all $i=1,...,s$.

Let now $\mu=\sum\limits_{i=1}^{s}m_{i}\delta_{x_{i}}\in P_{\omega}(X)$ and $m_{i}\leq \frac{1}{2}$ for all $i=1,...,s$. Fix a set $\{x_{1},...,x_{k}\}$. Then points of the view $\frac{1}{2}\delta_{x_{i}}+\frac{1}{2}\delta_{x_{j}}$, $i,j\in\{1,...,k\}$ are vertices of convex set formed by affine combinations $\sum_{i=1}^{s}m_{i}\delta_{x_{i}}$  Dirac measures $\delta_{x_{i}}$, and numbers $m_{i}$,  $0\leq m_{i}\leq\frac{1}{2}$,  $i=1,...,k$,  $\sum_{i=1}^{s}m_{i}=1$. Hence $\mu\in P_{\frac{\omega}{2}}(X)$. Claim 3 is proved.

Return now to proof Lemma 1. Thus a set  $P_{\omega}(X)\backslash P_{\frac{\omega}{2}}(X)$ is open in $P_{\omega}(X)$ and $P_{f}(X)\subset P_{\omega}(X)\backslash P_{\frac{\omega}{2}}(X)$. Let $\sum\limits_{i=1}^{m}\alpha_{i}\delta_{x_{i}}\in P_{\frac{\omega}{2}}(X)$. Then $\alpha_{i_{0}}>\frac{1}{2}$ for some (a unique) $i_{0}\in \{1,...,m\}$. We will construct $r:P_{\omega}(X)\backslash P_{\frac{\omega}{2}}(X)\rightarrow P_{f}(X)$ as follows
\[
r\left(\sum_{i=1}^{n}\alpha_{i}\delta_{x_{i}}\right) =
\begin{cases}
\frac{n}{n+1}\cdot\delta_{x_{i_{0}}}+\sum_{\substack{i=1,\\ i\neq i_{0}}}^{n}\frac{\alpha_{i}}{(n+1)(1-\alpha_{i_{0}})}\delta_{x_{i}}, & \text{at $\alpha_{i_{0}}<\frac{n}{n+1}$;} \\
\sum_{i=1}^{n}\alpha_{i}\delta_{x_{i}}, & \text{at $\alpha_{i_{0}}\geqslant \frac{n}{n+1}$.}
\end{cases}
\]
It is easy to see that the map $r$ is retraction. Lemma 1 is proved.

{\bf Theorem 3.} {\it Let $X$ be an $A(N)R$-compact. Then $P_{f}(X)$ is $ANR$-compact.}

{\bf Proof.} Let $X$ be a neighbourhood retract of some compact $Y$, $U$ be an open set in $Y$, such that $U\supset X$ and exist retraction $r:U\rightarrow X$.

Consider open set $\left\langle U;\frac{1}{2}\right\rangle=\left\{\mu\in P_{f}(Y):\mu(U)>\frac{1}{2}\right\}$ in $P_{f}(Y)$. It is evident that $P_{f}(X)\subset \left\langle U;\frac{1}{2}\right\rangle$. Let $\nu\in\left\langle U;\frac{1}{2}\right\rangle\subset P_{f}(Y)$,  $\nu=\sum_{i=1}^{n}n_{i}\delta_{y_{i}}$,  $\sum_{i=1}^{n}n_{i}=1$,  $n_{i}\geq 0$,  $i=1,2,...,n$,  $n_{i_{0}}\geq\frac{n}{n+1}$. Then $r_{f}^{x}(\nu)=\delta_{y_{i_{0}}}$. It is easy to check that $y_{i_{0}}\in U$, and hence $\delta_{y_{i_{0}}}\in\left\langle U;\frac{1}{2}\right\rangle$. Put
$$r_{U}^{Y}(\nu)=\left(n_{i_{0}}+\sum_{y_{i}\in Y\backslash U}n_{i}\right)\delta_{y_{i_{0}}}+\sum_{y_{i}\in U}n_{i}\delta_{y_{i}}.$$
It is obvious that $r_{U}^{Y}(\nu)\in P_{f}(Y)$  and  $supp r_{U}^{Y}(\nu)\subset U$. Besides, $r_{U}^{Y}(\nu)=\nu$ for every measure $\nu\in P_{f}(Y)$ such that  $supp\nu\subset U$. It is easy to establish that map $r_{U}^{Y}(\nu):\left\langle U;\frac{1}{2}\right\rangle\rightarrow \left\langle U;\frac{1}{2}\right\rangle$ is continuous. Further, put
$$R\left(r_{U}^{Y}(\nu)\right)=\left(n_{i_{0}}+\sum_{y_{i}\in Y\backslash U}n_{i}\right)\delta_{r\left(y_{i_{0}}\right)}+\sum_{y_{i}\in U}n_{i}\delta_{r\left(y_{i}\right)}.$$
By construction inclusion $R\left(r_{U}^{Y}(\nu)\right)\in P_{f}(X)$ holds. The map $R\left\langle U;\frac{1}{2}\right\rangle\rightarrow P_{f}(X)$ is defined correctly. Since the retraction $r:U\rightarrow X$ is continuous, then the map $R$ is also continuous. Obviously, that $R\left(r_{U}^{Y}(\nu)\right)=\nu$ for every measures $\nu\in P_{f}(X)$. Thus, $R\circ r_{U}^{Y}:\left\langle U;\frac{1}{2}\right\rangle\rightarrow P_{f}(X)$  is required retraction. So, the set $P_{f}(X)$ is neighbourhood retract of the compact $P_{f}(Y)$.

Now using Lemma 1 and Theorem (3.1) (see [7]) completes the proof of Theorem 3.

On functorial language Theorem 3 looks as so:

{\bf Corollary 1.} {\it Functor $P_{f}$ preserves $ANR$-compacts.}

Theorem 3 and Proposition 1 imply the main result of the paper.

{\bf Corollary 2.} {\it Let $X$ be a compact. $P_{f}(X)\in ARN$ if and only if $X\in ARN$.}

Further, Theorem 3 implies the following statements.

{\bf Corollary 3.} {\it Functor $P_{f}$  preserves property of a compacts to be $Q$-manifold or Hilbert cube.}

{\bf Corollary 4.} {\it Functor $P_{f}$ preserves property of fibres of maps to be $ANR$-compact, compact  $Q$-manifold and Hilbert cube (finite sum of Hilbert cube).}

\begin{center}
\textsl{References}
\end{center}

{1.} \v{S}\v{c}epin, E. V. Functors and uncountable degrees of compacta. (Russian) Uspekhi Mat. Nauk 36 (1981), no. 3(219), 3–62, 255.

{2.} Djuraev, T. F. Some main properties of the functor $P_f$. (Russian) Vestnik Moscow State University. Series  Math.Mech. 1989. No 6. P. 29-33.

{3.} Djuraev, T. F. Space of all probability measures with finite support is homeomorphic to infinite dimensional linear space. (Russian) General topology. Spaces and mappings. Moscow. MSU. 1989. P. 66-70.

{4.} Djuraev, T. F. On functor $P$ of probability measures. (Russian) Vestnik Moscow State University. Series  Math.Mech. 1990 No 1. P. 26-30.

{5.} Djuraev, T. F. Some geometrical properties of subfunctors of the functor $P$ of probability measures. (Russian). Moscow. MSU. 1989. 60 p. The manuscript is deposited at VINITI of Academy of Sciences of the USSR bibliography 58 names. 1989, July 5. No 4471 -- В 89.

{6.} Borsuk, K. Theory of shape. Monographie Mat., Vol. 59. PWN (Polish scientific publishers). Warszawa. 1975.

{7.} Borsuk, K. Theory of retracts. Monographie Mat., Vol. 44. PWN (Polish scientific publishers). Warszawa. 1967.

{8.} Varadarajan, V. S. Measures on topological spaces. (Russian) Mat. Sb. (N.S.) 55 (97) 1961 35-100.

{9.} Chapman, T. A. Lectures on Hilbert Cube Manifolds. CBMS Regional Conference Series in Mathematics
Volume: 28; 1976; 131 pp.

\end{document}